# Fractional Order Load-Frequency Control of Interconnected Power Systems Using Chaotic Multi-objective Optimization

Indranil Pan[a,b] and Saptarshi Das[b,c,*]

a) *Centre for Energy Studies, Indian Institute of Technology Delhi, Hauz Khas, New Delhi 110016, India.*

b) *Department of Power Engineering, Jadavpur University, Salt Lake Campus, LB-8, Sector 3, Kolkata-700098, India.*

c) *School of Electronics and Computer Science, University of Southampton, Southampton SO17 1BJ, United Kingdom.*

**Authors' Emails:**

indranil.jj@student.iitd.ac.in, indranil@pe.jusl.ac.in (I. Pan)

saptarshi@pe.jusl.ac.in, s.das@soton.ac.uk (S. Das*)

Phone: +44-7448572598

**Abstract:**

Fractional order proportional-integral-derivative (FOPID) controllers are designed for load frequency control (LFC) of two interconnected power systems. Conflicting time domain design objectives are considered in a multi objective optimization (MOO) based design framework to design the gains and the fractional differ-integral orders of the FOPID controllers in the two areas. Here, we explore the effect of augmenting two different chaotic maps along with the uniform random number generator (RNG) in the popular MOO algorithm – the Non-dominated Sorting Genetic Algorithm-II (NSGA-II). Different measures of quality for MOO e.g. hypervolume indicator, moment of inertia based diversity metric, total Pareto spread, spacing metric are adopted to select the best set of controller parameters from multiple runs of all the NSGA-II variants (i.e. nominal and chaotic versions). The chaotic versions of the NSGA-II algorithm are compared with the standard NSGA-II in terms of solution quality and computational time. In addition, the Pareto optimal fronts showing the trade-off between the two conflicting time domain design objectives are compared to show the advantage of using the FOPID controller over that with simple PID controller. The nature of fast/slow and high/low noise amplification effects of the FOPID structure or the four quadrant operation in the two inter-connected areas of the power system is also explored. A fuzzy logic based method has been adopted next to select the best compromise solution from the best Pareto fronts corresponding to each MOO comparison criteria. The time domain system responses are shown for the fuzzy best compromise solutions under nominal operating conditions. Comparative analysis on the merits and de-merits of each controller structure is reported then. A robustness analysis is also done for the PID and the FOPID controllers.

1. Introduction

Large scale power system networks comprise of several interconnected subsystems representing particular geographical areas. Each of these subsystems has their own generation capability and has variable load demand. These sub-systems are connected by tie lines which control the flow of power between the different areas [1]. A sudden load demand in a certain area results in a drop in system frequency which is detrimental for connected electrical loads. To ensure proper power quality, the load frequency controllers in the interconnected power system regulate the flow of power among the different areas through the tie lines and balance the load and the drop in frequency. The LFCs balance the mismatch between the frequencies of the interconnected areas and schedule the flow of power through the tie lies, helping the interconnected power system to overcome the aberrations introduced due to varying load demand, generation outage etc. Recently the LFCs are gaining more importance due to the integration of renewables in the grid which have an inherent stochastic characteristic due to the vagaries of nature unlike those of the base load thermal power plants [2], [3]. Thus proper design and operation of the LFCs are very important for the stable and reliable operation of large scale power systems. Control of interconnected power system considering various aspects has been a topic of intense research in the recent past. Different type of generating units and their effects have been studied e.g. thermal with reheat [4], generation rate constraint (GRC) [5], reheat and battery energy storage both the areas [6], hydro turbine and hydro-governor in both the areas [7], thermal with reheat turbine along with hydro and gas turbine plants in both the areas [8], etc.

Traditionally a proportional-integral (PI) or a proportional-integral-derivative (PID) controller is used for the LFCs [9] and a variety of different methods exist for proper tuning of the controller parameters. Many robust control design techniques have been applied to the LFC problem so that the designed controller is able to handle uncertainties of the system. Some variants of robust designs include an adaptive output feedback based robust control [10], adaptive robust control [11], decentralized robust control using iterative linear matrix inequalities (LMIs) [12], decentralized control [13] etc. Optimal control designs using Linear Quadratic Regulator (LQR) technique has been reported in [14]. Several other popular control philosophies like the model predictive control (MPC) [15], sliding mode control [16], and singular value decomposition (SVD) [17], etc. have also been applied in decentralized LFC of multi-area power systems.

Many computational intelligence based techniques have been employed in the design of LFCs as well. Global optimization techniques using evolutionary and swarm intelligence has been used to tune the PID controller parameters in various literatures. Genetic algorithm (GA) has been used in the design of LFCs in [7]. A variable structure controller has been designed for LFC's using GA in [18]. Fuzzy logic based gain scheduling has been done in [19] to obtain improved control strategy for LFC. The application of neural networks in the problem of load frequency control has been investigated in [20], [21]. There are other literatures which employ robust control techniques like $H_\infty$ loop shaping [22], $\mu$-synthesis [23] and LMI approaches [12], [22] with intelligent genetic algorithms to obtain robust controllers. A detailed review of the existing design methodologies in LFC has been documented in [24]-[26]. Other intelligent algorithms like Bacterial foraging algorithm [27], fuzzy logic [28], recurrent fuzzy neural network [29] have also been applied in LFC of multi-area inter-connected power systems.

Fractional order controllers have been gaining popularity in recent years due to added capability to handle control design specifications [30], [31]. Fractional order PID controllers have been applied in a wide variety of control systems and have generally proven to be better than their integer order (IO) counterparts [32]. Recently fractional order controllers have been applied to power systems and favourable results have been obtained. In [33] a fractional order controller has been designed for an automatic voltage regulator (AVR) with particle swarm optimization (PSO) algorithm to show that the FO controllers have more robustness to tackle uncertainties than the conventional IO-PID controller. Alomoush [34] has applied fractional order controllers for LFC of a two area interconnected power system and an automatic generation control for an isolated single area power system. It has been shown in [34] that the fractional order PID controller has more flexibility in design and can adjust the system dynamics better than the IO-PID controller. FOPID controllers are also shown to be robust and competitive to IO-PID controllers [35–38]. In these previous literatures, only a single objective intelligent optimization has been employed to design the control system. However it is well known that there exists multiple trade-offs among different design specifications in control [39] and similar system design [40]. It is not possible to simultaneously minimize all design objectives using a particular control structure and different controller structure may yield different trade-offs depending on the choice of the conflicting control objectives [41], [42]. Thus there is a requirement of multi-objective approach for addressing different conflicting objectives in the control system design [42]. In [42–44], time and frequency domain multi-objective formulation have been used to study the design trade-offs among various conflicting FOPID design objectives in an automatic voltage regulator (AVR) system. It has been shown in [42–44] that sometimes the FOPID and at other times the PID controller performs better depending on the choice of the conflicting objective functions. This concept of MOO for FO controllers has been extended in this paper for two area LFC problem. In this paper, FOPID controllers in the two area LFC are designed using chaotic multi-objective NSGA-II algorithm. Set point tracking and low control signal are chosen as the two conflicting objectives for the MOO based tuning of FOPID/PID controllers and the performance improvement due to the FOPID compared to the PID is illustrated by numerical simulations.

The main highlights of the paper includes:

- Augmenting the NSGA-II algorithm with different chaotic maps (like Logistic and Henon maps) to obtain better Pareto optimal solutions.

- Using Pareto metrics like hyper-volume indicator, spacing metric, Pareto spread and diversity metric [45]-[47], to assess the performance of the chaotic MOO algorithms.

- Use of FOPID controller to obtain better system performance for the two area LFC.

- Demonstration of conflicting time domain trade-offs in the controller performances for the FOPID and the PID controller and use of a fuzzy based mechanism for selecting the best compromise solution.

- Robustness study of FOPID as LFC over that with PID, under system parametric uncertainty and random change in load patterns.

The rest of the paper is organised as follows. Section 2 gives a brief description of the two area LFC problem in the interconnected power system. Section 3 introduces the basics of fractional calculus, FOPID controller and its flexibility over PID. Section 4 outlines the conflicting time domain criteria based MOO for the design of the LFC system. Section 5 introduces the chaotic versions of the multi-objective NSGA-II algorithm. Different MOO measure based selection of the best optimizer and the controller are enunciated in Section 6, along with the respective time domain responses. In section 7, the effect of uncertainty in the synchronizing coefficient and the effect of randomly changing load patterns in both areas are explored. The paper ends with the conclusions in Section 8 followed by the references.

## 2. Load frequency control of interconnected two area power system

The main functionalities of the LFC are

a) To keep the operating power system frequency within specified tolerance limits to ensure power quality
b) To ensure proper load sharing between the generators of the interconnected system
c) To honour the pre-specified load exchange constraints by controlling the power flow between the interconnected areas.

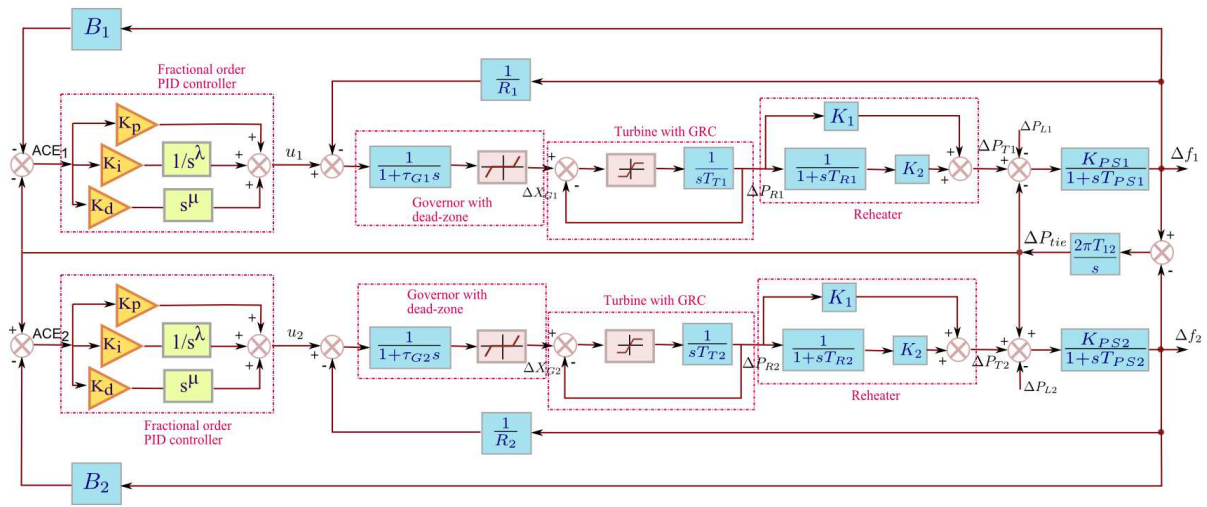

Figure 1: Block diagram representation of a two area AGC system with secondary LFC loop

Table 1: Data for the two area system considered in the present simulation

| Area | $K_{PS}$ (Hz/pu) | $T_{PS}$ (s) | R (Hz/pu MW) | B (pu MW/Hz) | $T_G$ (s) | $T_T$ (s) | $T_R$ (s) | $K_1$ | $K_2$ | $T_{12}$ |
|---|---|---|---|---|---|---|---|---|---|---|
| Area 1 | 120 | 20 | 2.4 | 0.425 | 0.1 | 0.3 | 10 | 0.5 | 0.5 | 0.0707 |
| Area 2 | 120 | 20 | 2.4 | 0.425 | 0.1 | 0.3 | 10 | 0.5 | 0.5 | |

A schematic diagram for the two area LFC is depicted in Figure 1. The parameters of the various units of the system are shown in Table 1 [23]. The area control error (ACE) in each area is a function of the frequency deviation ($\Delta f$) and the inter area tie-line power flow ($\Delta P_{tie}$). The ACE is fed as an input to the FOPID controller where it calculates the appropriate control signal to be applied. This control signal is fed into the amplifier and then the actuator to produce an appropriate change in the mechanical torque of the turbine (prime mover). This produces a change in the active power output of the generator to compensate the power flow in the system and thus $\Delta f$ and $\Delta P_{tie}$ are kept within desired limits. The tasks of the PID/FOPID controller in each area are to ensure faster damping of individual ACEs and also damping the inter-area power oscillation $\Delta P_{tie}$.

Each area includes steam turbine, governor, reheater stages along with GRC nonlinearity in the turbine and dead-band in the governor. The power systems (PS) are represented by first order transfer functions. The tie-line power is also affected by the choice of synchronizing coefficient ($T_{12}$). There are different configurations like primary and secondary LFC as described in [34]. In the primary LFC control, load change in one area is not corrected by a controller in that area or in other areas. For the secondary LFC loop, the tie-line is represented by the accelerating power coefficient ($P_s = 2\pi T_{12}$). Any change in the demand-load ($\Delta P_L$) will result in deviation of frequency in both the areas and the tie-line power flow. The ACE for the $i^{th}$ area can be expressed as (1).

$$ACE_i = \sum_{j=1}^{M} \Delta P_{ij} + B_i \Delta f_i \tag{1}$$

where $\Delta P_{ij}$ is the deviation in tie-line power flow from its scheduled values between the $i^{th}$ area and the $j^{th}$ area, $\Delta f_i$ is the frequency aberration in the $i^{th}$ area and M is the number of areas connected to area $i$. The frequency bias factor ($B_i$) can be expressed as a combination of the speed regulation ($R_i$) and the damping coefficient ($D_i$) and is given by (2).

$$B_i = (1/R_i) + D_i \tag{2}$$

In the present model, significant effect of nonlinearity is studied in the form of dead-zone in the governors and GRC in the turbines. The GRC keeps the rate of change of power within a specified limit of $\delta = \pm 0.005$ which is implemented by replacing the linear model of the turbine ($\Delta P_R/\Delta X_G$) by the nonlinear one as shown in Figure 1. The governor dead-band affects the speed control under disturbances and has been chosen as 0.06% in each area. For the present simulation study, both the areas are subjected to step-load disturbance of $P_{L1} = 0.02\,pu$ and $P_{L2} = 0.008\,pu$ respectively. In contemporary literature there were several studies on the dynamics of nonlinearities like GRC [5], [29], [37], [36] and dead-zone [15] along with reheat turbine [48], but their MOO based FO control design has not been investigated yet, which is the main motivation of the present paper.

### 3. Fractional calculus and Fractional order PI$^\lambda$D$^\mu$ (FOPID) controller

The generalized fractional differentiation and integration has mainly three definitions, the Grunwald-Letnikov definition, Riemann-Liouville definition and Caputo definition. The Grunwald-Letnikov formula is basically an extension of the backward finite difference formula for successive differentiation. This formula is widely used for the numerical solution of fractional differentiation or integration of a function. The Riemann-Liouville definition is an extension of *n*-fold successive integration and is widely used for analytically finding fractional differ-integrals. In the FO systems and control related literatures, mostly the Caputo's definition of fractional differ-integration is referred. This typical definition of fractional derivative is generally used to derive fractional order transfer function models from fractional order ordinary differential equations with zero initial conditions. According to Caputo's definition, the $\alpha^{th}$ order derivative of a function $f(t)$ with respect to time is given by (3) which is used in the present paper for realizing the fractional integro-differential operators of the FOPID controllers.

$$^C_0 D^\alpha_t f(t) = \frac{1}{\Gamma(m-\alpha)} \int_0^t \frac{D^m f(\tau)}{(t-\tau)^{\alpha+1-m}} d\tau, \alpha \in \mathbb{R}_+, m \in \mathbb{Z}_+, m-1 < \alpha < m. \qquad (3)$$

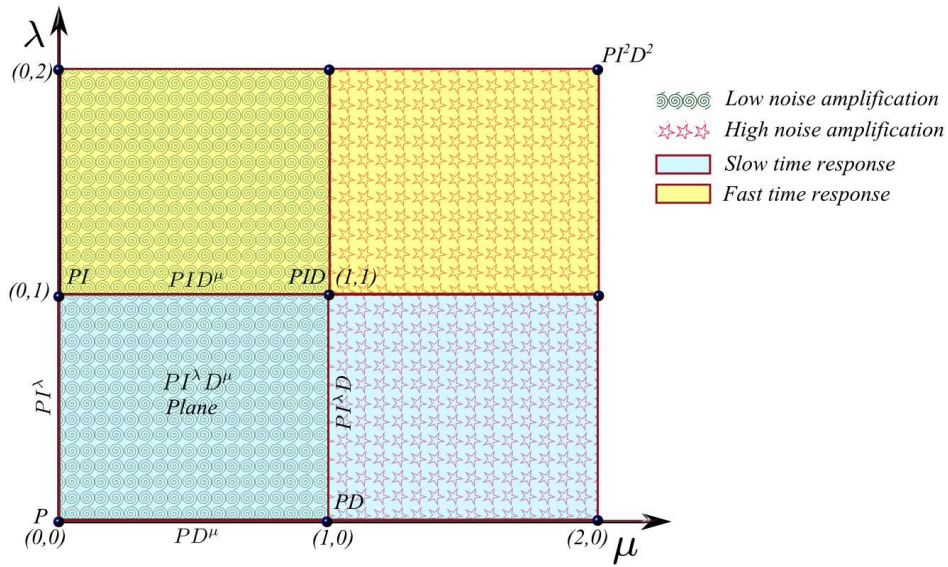

**Figure 2: Four quadrant operation of FOPID with different choice of speed level and noise amplification, compared to the conventional integer order PID controller**

The FOPID controller is an extension of the IO-PID controller with non-integer choice of integro-differential orders along with the conventional PID controller gains. The transfer function representation of a FOPID controller is given in (4).

$$C(s) = K_p + \left(K_i / s^\lambda\right) + K_d s^\mu \qquad (4)$$

This typical controller structure has five independent tuning knobs i.e. the three controller gains $\{K_p, K_i, K_d\}$ and two fractional order integro-differential operators $\{\lambda, \mu\}$. For $\lambda = 1$ and $\mu = 1$, the controller structure (4) reduces to the classical PID controller in parallel structure. Figure 2 shows the schematic representation of the FOPID controller in the integro-differential (λ-μ) plane and its relation with the conventional integer order P, PI, PD and PID controllers. Several other special cases of the FO controller structure can also be defined in the λ-μ plane like $PI^\lambda$, $PD^\mu$, $PI^\lambda D$, $PID^\mu$ etc. It has been shown in [31] that increasing the integral order above $\lambda > 1$ amplifies the low frequency components thus making the closed-loop system oscillatory. Similarly, an increase in the derivative order $\mu > 1$ increases the high frequency gain, thus amplifying the measurement noise and high frequency random components. The extension of 'points to plane' has led to the concept of FOPID or $PI^\lambda D^\mu$ controllers as shown in Figure 2 which can be further refined to define regions of low/high noise amplification and slow/fast time response depending on the range of values for μ and λ, below/above one. Figure 2 shows how different trade-offs could be achieved by selecting the range of operation in the four quadrants of the λ-μ plane (instead of only one quadrant with $\{\lambda, \mu\} < 1$), for different application of FOPID controller. For the present power system control application, the controller selection problem is divided in two parts – fast FOPID (λ>1 or 1st/2nd quadrant operation) and slow FOPID (λ<1 or 3rd/4th quadrant operation). Additionally, for choice of speed (λ), four other combinations of the controllers in the two areas are explored, depending on the derivative order μ>1 or μ<1. Within the present MOO framework, considering similar or different controller structure in both the areas, the best non-dominated solution generated by the four controller combinations is selected. The choice of the controller within the MOO framework has been shown in (5).

$$\begin{aligned}
&\text{Select Slow FOPID}\,(\lambda_1, \lambda_2 < 1), \text{in both the areas} \\
&\quad \left.\begin{array}{l}
\text{choice 1:} \quad \mu_1 < 1, \text{in area 1 and } \mu_2 < 1, \text{in area 2} \\
\text{choice 2:} \quad \mu_1 < 1, \text{in area 1 and } \mu_2 > 1, \text{in area 2} \\
\text{choice 3:} \quad \mu_1 > 1, \text{in area 1 and } \mu_2 < 1, \text{in area 2} \\
\text{choice 4:} \quad \mu_1 > 1, \text{in area 1 and } \mu_2 > 1, \text{in area 2}
\end{array}\right\} \text{selection by MOO} \\
&\text{Select Fast FOPID}\,(\lambda_1, \lambda_2 > 1), \text{in both the areas} \\
&\quad \left.\begin{array}{l}
\text{choice 1:} \quad \mu_1 < 1, \text{in area 1 and } \mu_2 < 1, \text{in area 2} \\
\text{choice 2:} \quad \mu_1 < 1, \text{in area 1 and } \mu_2 > 1, \text{in area 2} \\
\text{choice 3:} \quad \mu_1 > 1, \text{in area 1 and } \mu_2 < 1, \text{in area 2} \\
\text{choice 4:} \quad \mu_1 > 1, \text{in area 1 and } \mu_2 > 1, \text{in area 2}
\end{array}\right\} \text{selection by MOO}
\end{aligned} \qquad (5)$$

Therefore, within the slow and fast family of FOPID controllers, different levels of noise amplification in the respective areas are selected automatically by the MOO which results in non-dominated Pareto front.

Few recent research results show that band-limited implementation of FOPID controllers using higher order rational transfer function approximation of the integro-differential operators gives satisfactory performance in industrial automation. Here, the Oustaloup's recursive approximation has been used to implement the integro-differential operators in frequency domain is given by (6), which represents a higher order analog filter.

$$s^\alpha \simeq K \prod_{k=-N}^{N} \frac{s+\omega'_k}{s+\omega_k} \qquad (6)$$

where, the poles, zeros, and gain of the filter can be recursively evaluated as (7).

$$\omega_k = \omega_b \left(\omega_h/\omega_b\right)^{\frac{k+N+(1+\alpha)/2}{2N+1}}, \omega'_k = \omega_b \left(\omega_h/\omega_b\right)^{\frac{k+N+(1-\alpha)/2}{2N+1}}, K = \omega_h^\alpha \qquad (7)$$

Thus, the area control errors (ACEs) can be passed through the filter (6) and the output of the filter can be regarded as an approximation to the fractionally differentiated or integrated signal $D^\alpha f(t)$. These FO differentiated or integrated signals are weighted by the respective gains to form the final control signal which goes to the governor in Figure 1. In (6)-(7), $\alpha$ is the order of the differ-integration, $(2N+1)$ is the order of the filter and $(\omega_b, \omega_h)$ is the expected frequency fitting range. In the present study, 5$^{\text{th}}$ order Oustaloup's recursive approximation is implemented to approximate the integro-differential operators within the frequency band of $\omega \in \{10^{-2}, 10^2\}$ rad/sec for the constant phase elements (CPEs) of the FOPID controller.

## 4. Need of multi-objective optimization and conflicting time domain control objectives

It is well known that a single controller structure cannot give good results for all design specifications. For example, a fuzzy logic controller is good at coping with uncertainty in the loop, whereas a model predictive controller is good for tackling large time delays in process control. For specific applications, different controller structures would give a trade-off solution among conflicting design specifications. Hence for effective comparison of different controller structures it is essential to know the limits of performance of each of the individual controllers for conflicting design specifications. In [42]-[44], a similar approach has been taken to compare the efficacy of the FOPID controller vis-à-vis the IOPID one for several conflicting time and frequency domain objectives respectively. In the present case, the disturbance rejection and controller effort are considered as the conflicting objectives and expressed in (8)-(9).

$$J_1 = \sum_{i=1}^{M} ITSE_i = \sum_{i=1}^{M} \int_0^\infty t \cdot e_i^2(t) dt \qquad (8)$$

$$J_2 = \sum_{i=1}^{M} ISDCO_i = \sum_{i=1}^{M} \int_0^\infty \left(\Delta u_i(t)\right)^2 dt \qquad (9)$$

where, ITSE represents the Integral of the Time multiplied Squared Error, ISDCO represents Integral of the Squared Deviation in Controller Output, $e_i(t)$ represents the error signal (ACE) in area $i$, $u_i(t)$ represents the control signal and $M$ represents the total number of areas.

The reason for considering these two as conflicting objective functions can be briefly explained as follows. To achieve a faster damping of ACEs and the grid frequency oscillation, it is essential that the controller gains should be higher. In other words, the controller must be able to exert much more control action on the power generating system so that the frequency oscillation settles within a short amount of time. However, the control signal should ideally be smaller to prevent actuator saturation and minimize the cost associated with sizing of a larger actuator. It can be inferred that both these objectives of small control signal, as well as faster damping of load-disturbances cannot be ideally obtained by a fixed set of parameters of the PID/FOPID controller. Thus there would exist a range of values for the tuning parameters of the controller $\{K_p, K_i, K_d, \lambda, \mu\}$, where the controller would show good load disturbance rejection at the cost of higher control signal and vice-versa. After a large number of trade-off solutions between the two chosen objectives are obtained, a compromise solution could be selected next for deciding the most optimal controller setting [42], [49].

## 5. Multi-objective controller design using chaotic maps

### 5.1. Chaotic multi-objective optimization

A generalized multi-objective optimization framework can be defined as follows:

Minimize $F(x) = (f_1(x), f_2(x), ..., f_m(x))$

$$F(x) = (f_1(x), f_2(x), ..., f_m(x)) \qquad (10)$$

such that $x \in \Omega$; where $\Omega$ is the decision space, $\mathbb{R}^m$ is the objective space, and $F: \Omega \to \mathbb{R}^m$ consists of $m$ real valued objective functions.

Let, $u = \{u_1, ..., u_m\}, v = \{v_1, ..., v_m\} \in \mathbb{R}^m$ be two vectors. $u$ is said to dominate $v$ if $u_i < v_i \ \forall i \in \{1, 2, ..., m\}$ and $u \neq v$. A point $x^* \in \Omega$ is called Pareto optimal if $\nexists \ x \,|\, x \in \Omega$ such that $F(x)$ dominates $F(x^*)$. The set of all Pareto optimal points, denoted by PS$_{MOO}$ is called the Pareto set. The set of all Pareto objective vectors, $PF = \{F(x) \in \mathbb{R}^m, x \in PS_{MOO}\}$, is called the Pareto Front or the set of non-dominated solutions. This implies that no other feasible objective vector exists which can improve one objective function without simultaneous worsening of some other objective function.

Multi-objective Evolutionary Algorithms (MOEAs) which use non-dominated sorting and sharing, have higher computational complexity. They use a non-elitist approach and require the specification of a sharing parameter. The non-dominated sorting genetic algorithm (NSGA-II) removes these problems and is able to find a better spread of solutions and better convergence near the actual Pareto optimal front [50], [51].

The NSGA-II algorithm converts different objectives into one fitness measure by composing distinct fronts which are sorted based on the principle of non-domination. In the process of fitness assignment, the solution set not dominated by any other solutions in the population is designated as the first front $F_1$ and the solutions are given the highest fitness value. These solutions are then excluded and the second non-dominated front from the remaining population $F_2$ is created and ascribed the second highest fitness. This method is iterated until all the solutions are assigned a fitness value. Crowding distances are the normalized distances between a solution vector and its closest neighbouring solution vectors in each of the fronts. All the constituent elements of the front are assigned crowding distances to be later used for niching. The selection is achieved in tournaments of size 2 according to the following logic.

a) If the solution vector lies on a lower front than its opponent, then it is selected.

b) If both the solution vectors are on the same front, then the solution with the highest crowding distance wins. This is done to retain the solution vectors in those regions of the front which are scarcely populated.

The optimization variables for the fractional order PID controller are the proportional-integral-derivative gains and the differ-integral orders, i.e. $\{K_p, K_i, K_d, \lambda, \mu\}$ for both the areas. For the IO-PID controller the optimization variables are the gains only i.e. $\{K_p, K_i, K_d\}$ for both the areas. In other words, the dimension of the decision space $\Omega$ is five for the FOPID controller and three for the PID controller. The population size is taken as $15 \times n_{var}$ and the algorithm is run until the cumulative change in fitness function value is less than the function tolerance of $10^{-6}$ or the maximum generations ($200 \times n_{var}$) are exceeded. Here, $n_{var}$ of parameters to be tuned by the MOO algorithm in both the areas and given by $n_{var}^{PID} = 6$ and $n_{var}^{FOPID} = 10$ for the two controller structures. The crossover fraction is taken as 0.8 and an intermediate crossover scheme is adopted. The mutation fraction is chosen as 0.2 and the Gaussian scheme is adopted. For choosing the parent vectors based on their scaled fitness values, the algorithm uses a tournament selection method with a tournament size of 2. The Pareto front population fraction is taken as 0.7. This parameter indicates the fraction of population that the solver tries to limit on the Pareto front [51]. For the MOO problem, the limits of $\{K_p, K_i, K_d\}$ are chosen as $[0,10]$ and the bounds of the differ-integral orders $\{\lambda, \mu\}$ are chosen within the range $[0,2]$ which is described in (5).

The uniformly distributed RNG is normally used for the crossover and mutation operations in the standard version of the NSGA-II algorithm [50][51]. However since the strength of evolutionary algorithms lies in the randomness of the crossover and mutation operators, many contemporary researchers have focussed on increasing the efficiency of these algorithms by incorporating different random behaviours through various techniques like stochastic resonance and noise [52], chaotic maps [53] etc. These different strategies can be classified as special cases of a broader principle of diversification which essentially entails a trade-off between exploration and exploitation [54]. In [55] it has been shown that the performance of single objective evolutionary algorithms increase if different types of chaotic maps are introduced instead of the uniform RNG for the crossover and mutation operations. It has also been demonstrated in [55] that, in general, using chaotic systems for the RNG in the crossover and mutation operations may yield better result than using RNG from a noisy sequence in terms of convergence and effectiveness of the algorithms in finding global minima. In [56]-[58] it has been shown that the multi-objective NSGA-II algorithm can be improved by using chaotic maps and gives better result than the original NSGA-II algorithm in terms of convergence and efficiency. This is due to the fact that the chaotic process introduces diversity in the solutions. In this paper, we adopt similar policy and use chaotic logistic map and chaotic Henon map to obtain comparable solutions and convergence with respect to the standard NSGA-II algorithm. The logistic map is one of the simplest discrete time dynamical systems exhibiting chaos. The equation for the logistic map is given in (11).

$$x_{k+1} = ax_{k+1}(1 - x_k) \qquad (11)$$

The Henon map is a discrete time dynamical system that exhibits chaotic behaviour. Given a point with co-ordinates $\{x_n, y_n\}$, the Henon map transforms it to a new point $\{x_{n+1}, y_{n+1}\}$ using the set of equations in (12).

$$\begin{aligned} x_{n+1} &= y_n + 1 - ax_n^2, \\ y_{n+1} &= bx_n \end{aligned} \qquad (12)$$

The map is chaotic for the parameters $a = 1.4$ and $b = 0.3$. It is actually a simplified model of the Poincare section of the Lorenz system. The output $y_{n+1}$ varies in different ranges depending on the initial seed $\{x_0, y_0\}$. Since the Henon map is used here as a random number generator, it must produce a random number in the range $[0,1]$. Hence the output of the map for different initial conditions have been scaled in the range $[0,1]$ as also done in [59].

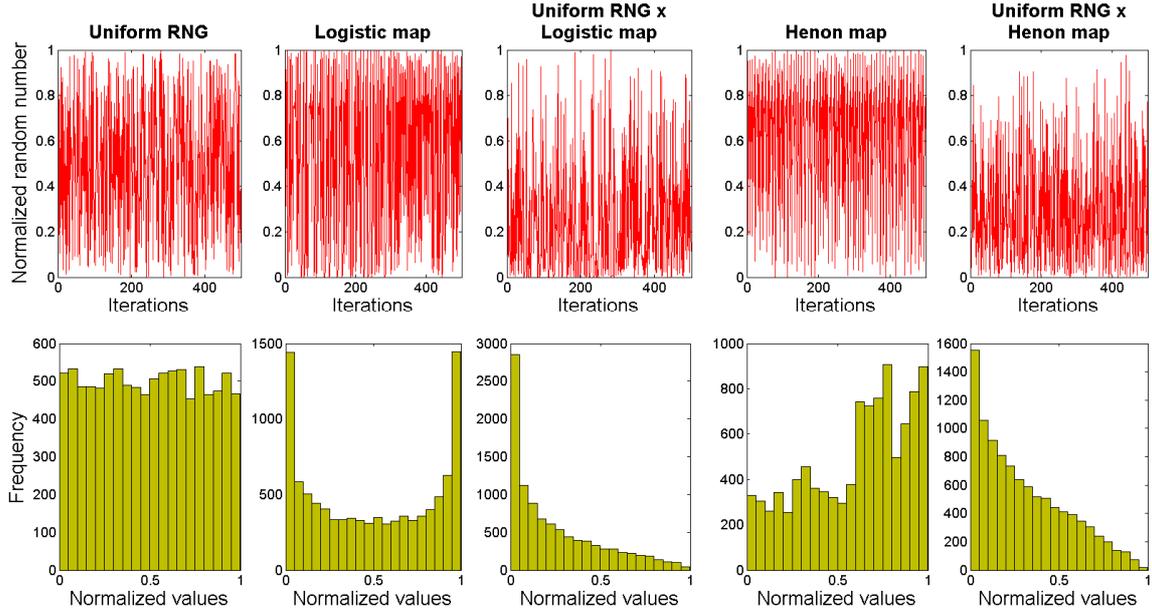

Figure 3: Random number generation by multiplying chaotic Henon map and Logistic map, normalized in the range [0,1] with uniform RNG. Top: first 500 samples of the sequences. Bottom: histogram of 10000 samples.

The initial seed of the Logistic map ($x_0$) in (11) has also been chosen randomly for the fixed parameter $a = 4$ exhibiting chaos, similar to that studied in [60]. For the implementation of the chaotic versions of the MOO algorithms, the outputs of the chaotic maps are scaled and multiplied with a uniform RNG which generates numbers in the range $[0,1]$. The first 500 samples of uniform RNG, the scaled output of the Henon map and the Logistic map and the outputs obtained by multiplying them with the uniform RNG are shown in Figure 3. It also shows the histogram of the respective cases with 10000 samples. It can be seen that even though the Logistic map and the Henon map have different distributions, the corresponding distributions obtained by multiplying them with a uniform RNG is highly skewed with lower values occurring more often than higher values. However, it is to be noted that even though both the histograms show similar skewed characteristics, it does not capture the time domain evolution like randomness and autocorrelation of the numbers. Therefore this generation mechanism with chaotic maps is different than simply drawing a random number from such skewed probability distributions each time.

## 5.2. Quantification of Pareto fronts and choosing the best compromise controller parameters

Due to the stochastic nature of the NSGA-II algorithm and also its two chaotic versions, the final Pareto fronts obtained at the end of each independent run would be slightly different with respect to one another. Therefore multiple runs are conducted to assess the convergence characteristics of the algorithms. There are several measures to compare multiple Pareto fronts apart from the non-domination criteria. In most realistic problems, two Pareto fronts under comparison could intersect with each other which indicates that in one region one set of solutions are more non-dominated in terms of one objective whereas for the other conflicting objective the other Pareto front would be more non-dominated. Especially in such cases of weak dominance, it could be difficult to judge the quality of the MOO solutions from a global perspective. In the case of strong Pareto dominance where one front dominates the other in all objectives, such problem is avoided due no intersection between multiple Pareto fronts. In order to avoid such case dependent strategy formulation, several generic metrics have been proposed in [45], which may indicate towards the quality of MOO solutions. In this paper, four different Pareto measures *viz.* minimum hypervolume indicator, maximum diversity metric, maximum Pareto spread, minimum spacing metric are explored. Since none of the measures can capture all necessary properties of the best possible Pareto front for the controller parameters [39], [61] a composite criteria is used to choose the best Pareto front.

The hypervolume indicator is the given by the total area/volume/hypervolume under the Pareto front with respect to the reference point as the origin. Therefore, a Pareto front closest to the origin will have the minimum hypervolume indicator showing a strong nondominance over other set of solutions [62]. The hypervolume indicator for two dimensions with reference to the origin is given by

$$HI = \sum_{i=1}^{N}(x_i - x_{i-1})(y_i - y_{i-1}) \qquad (13)$$

where, $N$ is the number of points on the Pareto front, $x, y$ are the two dimensions, $x_0, y_0$ correspond to the projections of the end points of the Pareto front on the two axes respectively. This has to be minimized to obtain a better non-dominated Pareto front.

The spacing metric measures the distance between the variance of neighbouring data-points on the Pareto-set [63], [64], which is given by (14).

$$SP = \sqrt{\frac{1}{S-1}\sum_{i=1}^{S}(\bar{d} - d_i)^2} \qquad (14)$$

where, the distance $d_i = \min_{j}\left(\sum_{k=1}^{m}|x_{ki} - x_{kj}|\right), \{i, j\} = 1, 2, \cdots, S$ and $\bar{d}$ is the mean value of $d_i$. Here, $S$ is the number of non-dominated solutions and $m$ is the number of objectives.

The total Pareto spread is given by the Euclidean distance between the extreme points of the Pareto front [65]. The total Pareto spread helps in understanding the span of the solutions along different conflicting objectives. For two dimensions, if $(x_1, y_1)$ and $(x_2, y_2)$ are the coordinates of the end points of the Pareto front, then the Pareto spread is given by

$$P_{spread} = \sqrt{(x_1 - x_2)^2 + (y_1 - y_2)^2} \qquad (15)$$

The moment of inertia based diversity metric is another popular measure to judge the quality of the Pareto fronts [66]. Let there be $S_{pt}$ number of points in the $m$-dimensional objective function space. Then the centroid for $i^{th}$ dimension is given by (16).

$$C_i = \left(\sum_{j=1}^{S_{pt}} x_{ij}\right) \bigg/ S_{pt}, \text{ for } i = 1, 2, \cdots, m \tag{16}$$

where, $x_{ij}$ denotes the $i^{th}$ dimension of the $j^{th}$ point. Then, the diversity metric can be calculated as (17).

$$I = \sum_{i=1}^{N} \sum_{j=1}^{S} \left(x_{ij} - C_i\right)^2 \tag{17}$$

Generally, the best, worst, mean and standard deviation of the diversity metric and other measure are compared for different MOO algorithms [46], [47]. In the present paper, we also report the box-plots of each of the four measures to ascertain the best Pareto front amongst 30 independent runs of each variant of the NSGA-II – nominal and chaotic. Additionally we report the comparison of the convergence times for different controller settings and optimizers.

After selecting the best Pareto front out 30 independent runs depending each of the four criteria i.e. minimum hypervolume indicator, maximum diversity metric, maximum Pareto spread, minimum spacing metric, a best compromise solution on the Pareto front has been obtained using a fuzzy based mechanism. In [43], [44] the median solution has been reported amongst all solutions on the Pareto front which has been improved here with a fuzzy based systematic choice of the best compromise solution [46], [47]. The designer may have imprecise goals for each objective function which can be encoded in the form of a fuzzy membership function $\mu_F$. Here, $\mu_{F_i}$ for each objective function $i$ is taken to be a strictly monotonic and decreasing continuous function expressed as (18).

$$\mu_{F_i} = \begin{cases} 1 & \text{if } F_i \leq F_i^{\min} \\ \left(F_i^{\max} - F\right)/\left(F_i^{\max} - F_i^{\min}\right) & \text{if } F_i^{\min} \leq F_i \leq F_i^{\max} \\ 0 & \text{if } F_i \geq F_i^{\max} \end{cases} \tag{18}$$

The value of $\mu_{F_i}$ represents the degree to which a particular solution has satisfied the objective $F_i$. The membership function lies between zero and one implying worst and best satisfaction of the objective respectively. The degree of satisfaction of each objectives by each solution can be represented as in (19).

$$\mu^k = \left(\sum_{i=1}^{m} \mu_i^k\right) \bigg/ \left(\sum_{j=1}^{S} \sum_{i=1}^{m} \mu_i^j\right) \tag{19}$$

where, $m$ is the number of objectives and $S$ is the number of solutions on the Pareto front. The best compromise solution on the Pareto front has been chosen in such a way for which (19) reaches its maximum.

## 6. Simulation and Results

### 6.1. Multi-objective criteria for the best controller selection

Similar to [23], the system in Figure 1 is simulated with a step input $\Delta P_{L1}$ of 0.02 pu in the first area and the $\Delta P_{L2}$ of 0.008 pu in the second area. The PID/FOPID controllers have been tuned using time domain performance indices in (8)-(9) under a MOO framework to handle these load inputs. Three different controller structures are considered – traditional PID controller, slow FOPID controller ($0 < \lambda < 1$) and the fast FOPID controller ($1 < \lambda < 2$). Each of these cases are run with three different optimization algorithms – the standard NSGA-II, the Logistic map adapted NSGA-II and the Henon map adapted NSGA-II. These 9 cases (3 controllers $\times$ 3 MOOs) are run for 30 times each and the statistics of the simulation time and Pareto metrics (like hypervolume indicator, diversity metric, Pareto spread and spacing metric) are calculated. The corresponding statistics are shown in the box-plots in Figure 4 to Figure 8. From Figure 4 it can be observed that the median time taken by the chaotic versions of the NSGA-II is higher than that of the normal NSGA-II for all the controllers. This is due to the additional computational time taken by the chaotic maps in the various crossover and mutation operators which employ RNG.

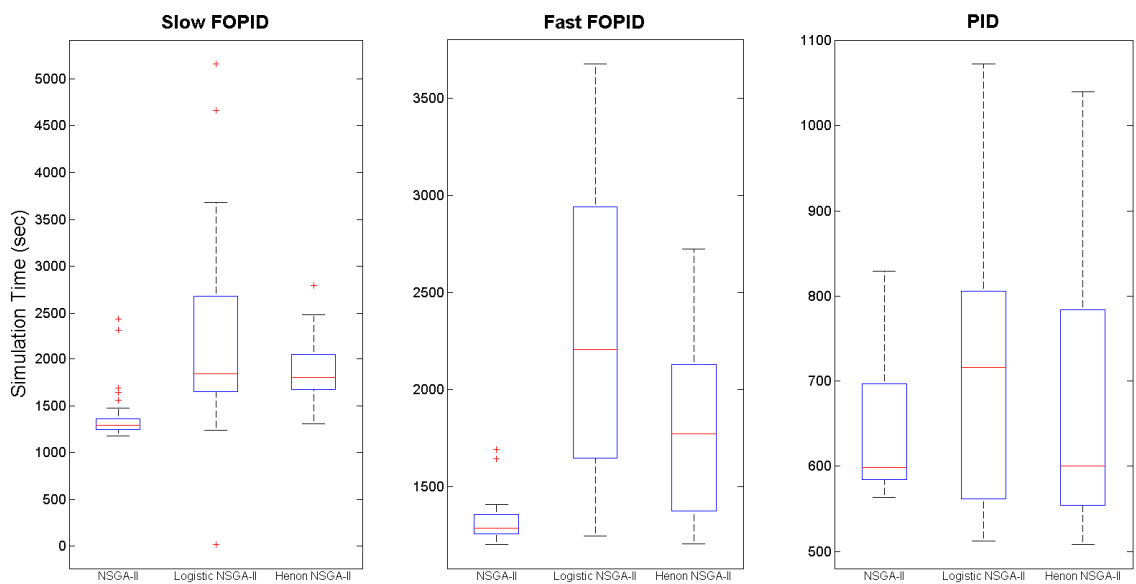

**Figure 4: Box-plot of the simulation times with different MOO algorithms and controller structures.**

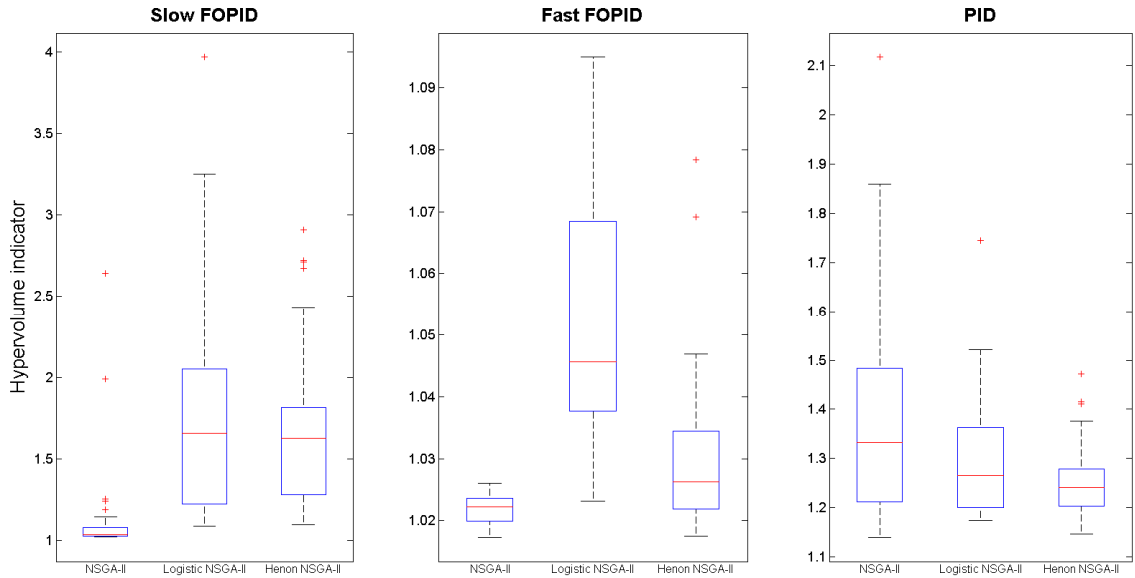

**Figure 5: Box-plot of the hypervolume indicator with different MOO algorithms and controller structures.**

Figure 5 shows the box plots of the hypervolume indicator for the three different algorithms and controller structures. From the minimum and median values, it is evident that the performance of the FOPID controller is better than the PID controller. The hypervolume indicator shows that the median values obtained using the traditional NSGA-II algorithm is better than the chaotic versions for the FOPID controllers but is worse for the PID structure.

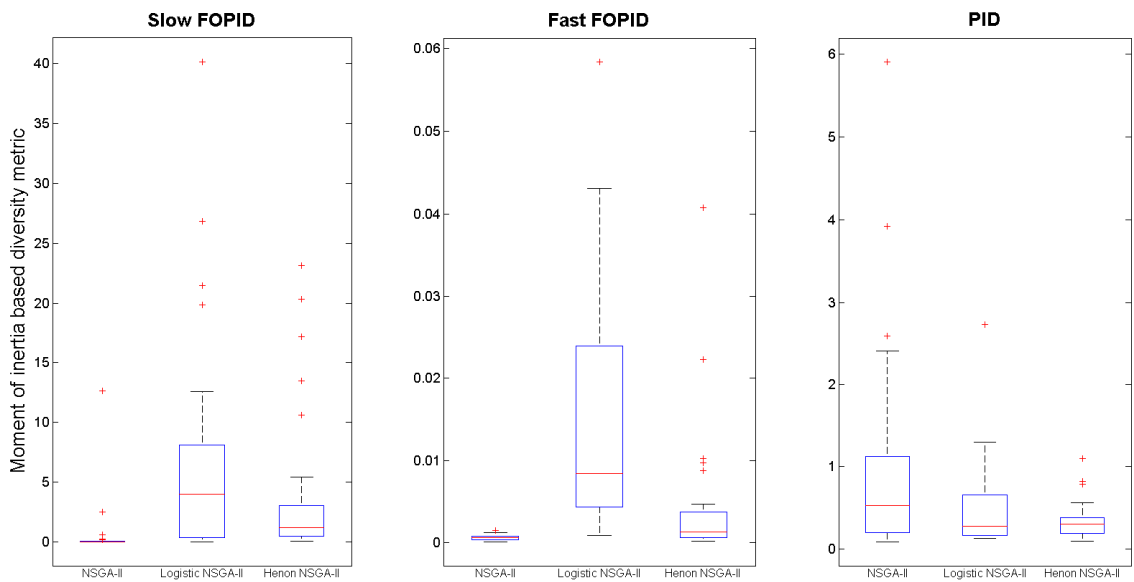

**Figure 6: Box-plot of the moment of inertia based diversity metric with different MOO algorithms and controller structures.**

Figure 6 shows the box plot of the moment of inertia based diversity metric for 30 runs of all the cases. A higher value of diversity metric indicates a better Pareto front. It can be observed that the Pareto fronts obtained for the FOPID controllers using the chaotic NSGA-II algorithms have a higher value of the diversity metric. But for PID controllers, the NSGA-II gives better results. It can also be seen that the slow FOPID controller has a higher range of the diversity metric than the PID controller. This is possibly due to the extra degrees of freedom of the FO integro-differential orders of the FOPID controller which has a larger search space and thus allows more diverse solutions. The fast FOPID has a very small value of diversity metric, indicating that most of the stable solutions belong to a small region of the search space. For all the other cases, the solutions are either unstable or dominated.

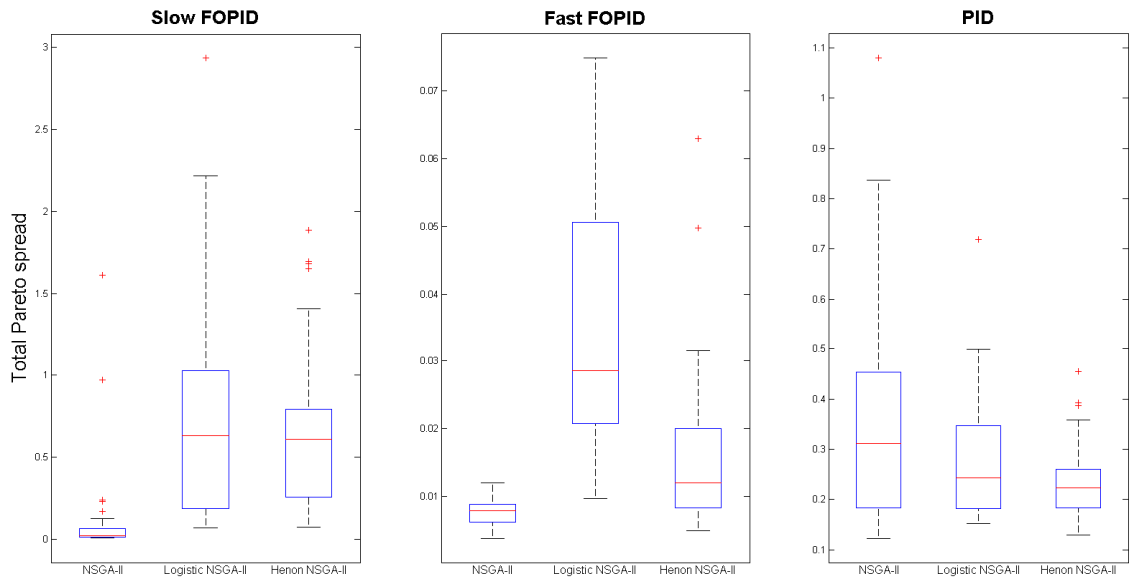

Figure 7: Box-plot of the total Pareto spread with different MOO algorithms and controller structures.

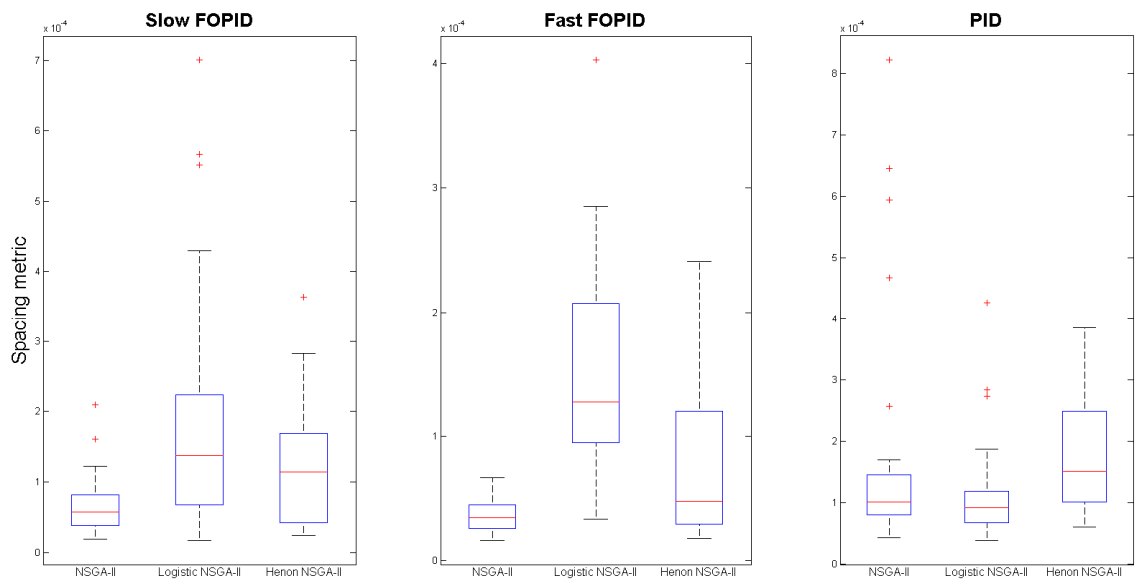

Figure 8: Box-plot of the spacing metric with different MOO algorithms and controller structures.

Figure 7 shows the box plots of the total Pareto spread for 30 runs of all the cases. It is observed that the chaotic versions of the NSGA-II are able to obtain a wider Pareto spread for the FOPID controllers, indicating a more diverse set of solutions. However, for the PID structure, the traditional NSGA-II gives a better Pareto spread. Figure 8 shows the box plots of the spacing metric for 30 runs of all the cases. The original NSGA-II is found to give a more uniform distribution of different solutions on the Pareto front.

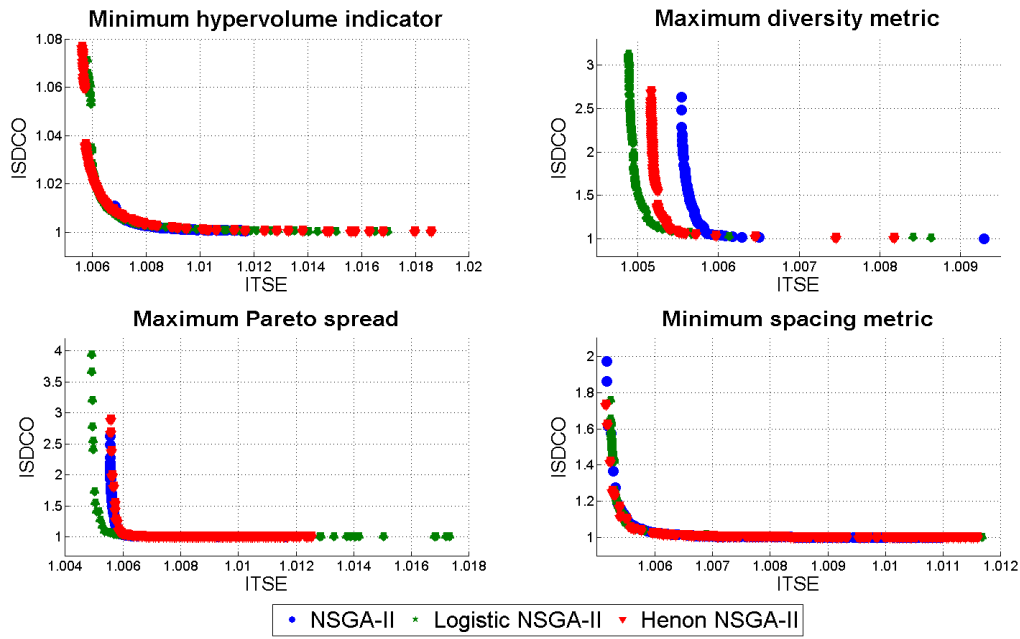

**Figure 9:** Comparison of the Pareto fronts for slow FOPID controller using different MOO algorithms and selection criteria.

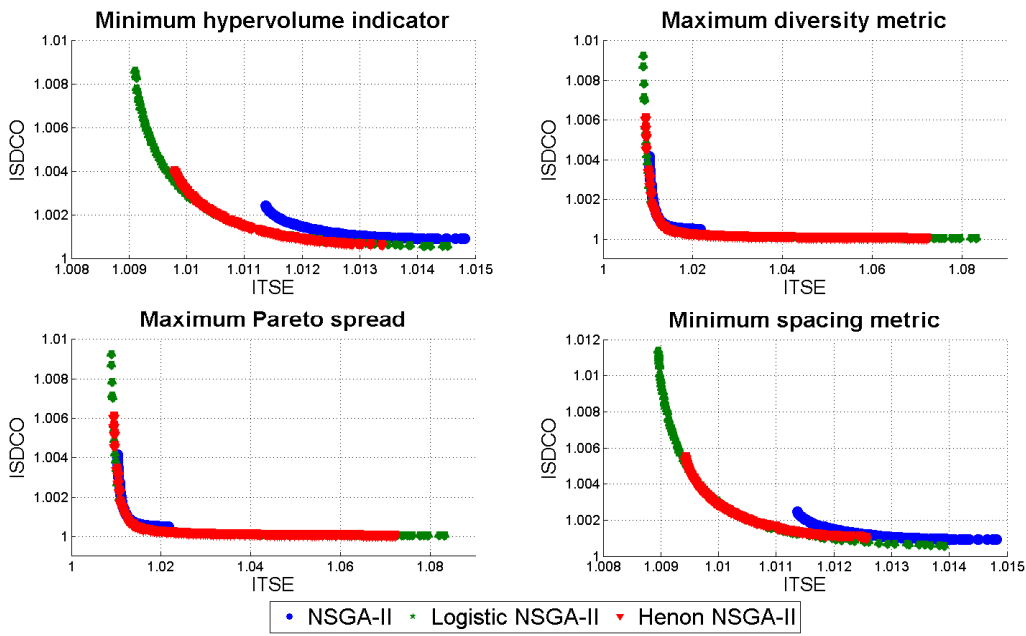

**Figure 10:** Comparison of the Pareto fronts for fast FOPID controller using different MOO algorithms and selection criteria.

It is therefore clear that none of the algorithms are better than their counterparts in all the metrics individually. Also each metric represents different characteristics of the Pareto front and finding the best Pareto front must leverage on some of these criteria taken together. The diversity metric, Pareto spread and spacing metric reflect the distribution of the solutions on the Pareto front. The hypervolume criterion which indicates the non-domination of the different fronts, directly affects the quality of the obtained solutions and a better non-domination implies a better control system performance. Hence non-domination is one of the most important metric among these. Therefore for comparison, the best Pareto front which is obtained by each of the different Pareto metrics is found out of each of the 30 runs.

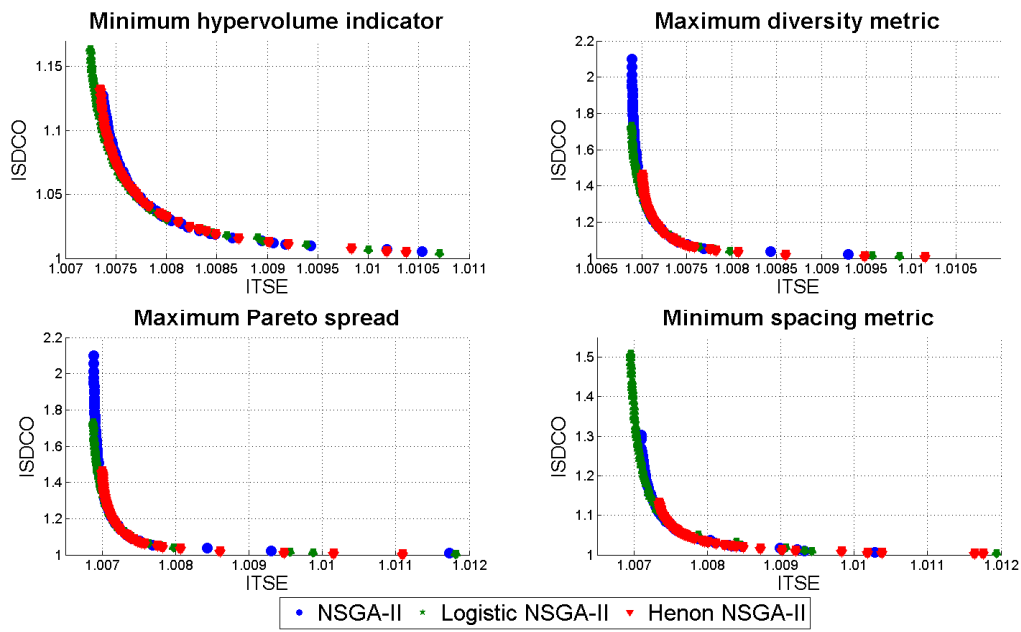

Figure 11: Comparison of the Pareto fronts for PID controller using different MOO algorithms and selection criteria.

It can be observed that for all the different Pareto metrics (hypervolume indicator, diversity metric, Pareto spread and spacing metric) the chaotic version of the NSGA-II (either logistic map assisted or Henon map assisted) gives the best non-dominated Pareto front. This can be verified from the superimposed Pareto fronts in Figure 9-Figure 11 for the three controller structures respectively. In most of the cases the logistic map assisted NSGA-II gives better performance than the others. In Figure 9-Figure 11, the best Pareto front (out of 30 runs) according to each of the four Pareto measures are shown. According to each criteria like the minimum hypervolume indicator, maximum diversity metric, maximum Pareto spread and minimum spacing metric, the chaotic NSGA-II versions gives wider and non-dominated Pareto spreads over that with the standard NSGA-II.

Next, for the sake of comparison amongst the best controller structures, the non-dominated Pareto fronts according to each of the Pareto metrics is identified from Figure 9-Figure 11 and are superimposed in Figure 12. It can be observed that under the nominal operating condition of the power system, the slow FOPID controller structure is more non-dominated and is therefore capable of resulting in better control performance. Although in all the four cases in Figure 12, the slow FOPID gives the non-dominant Pareto front but depending on different criteria the best compromise solution may perform better/worse. Therefore, we here report all the best compromise solution of the three controller structures for each of the Pareto fronts in Figure 12 (4 metrics × 3 controllers = 12 solutions in total). The corresponding algorithm which gives the best non-dominated Pareto front and along with the optimum controller parameters and objective functions have been reported in Table 2.

Table 2: Best compromise solutions of the FOPID and PID controller based on different Pareto metrics

| Controller | Criterion | Best Nondominated Algorithm | case | ITSE | ISDCO | $K_{p1}$ | $K_{i1}$ | $K_{d1}$ | $\lambda_1$ | $\mu_1$ | $K_{p2}$ | $K_{i2}$ | $K_{d2}$ | $\lambda_2$ | $\mu_2$ |
|---|---|---|---|---|---|---|---|---|---|---|---|---|---|---|---|
| Slow FOPID | Minimum hypervolume indicator | Henon NSGA-II | 1 | 1.01380 | 1.00040 | 0.090 | 0.297 | 0.036 | 0.869 | 0.208 | 0.036 | 0.237 | 0.036 | 0.533 | 0.585 |
| | Maximum diversity metric | Logistic NSGA-II | 2 | 1.01645 | 1.00053 | 0.215 | 0.230 | 0.026 | 0.573 | 0.724 | 0.160 | 0.121 | 0.068 | 0.354 | 0.612 |
| | Maximum Pareto spread | Logistic NSGA-II | 3 | 1.01145 | 1.00053 | 0.036 | 0.388 | 0.048 | 0.928 | 0.398 | 0.000 | 0.237 | 0.135 | 0.823 | 0.660 |
| | Minimum spacing metric | Henon NSGA-II | 4 | 1.01126 | 1.00068 | 0.010 | 0.484 | 0.083 | 0.570 | 0.544 | 0.090 | 0.165 | 0.140 | 0.523 | 0.715 |
| Fast FOPID | Minimum hypervolume indicator | Logistic NSGA-II | 5 | 1.01236 | 1.00082 | 0.129 | 0.397 | 0.107 | 1.014 | 0.553 | 0.244 | 0.173 | 0.139 | 1.155 | 0.745 |
| | Maximum diversity metric | Logistic NSGA-II | 6 | 1.04872 | 1.00007 | 0.008 | 0.100 | 0.037 | 1.055 | 0.203 | 0.034 | 0.036 | 0.014 | 1.222 | 0.699 |
| | Maximum Pareto spread | Logistic NSGA-II | 7 | 1.04872 | 1.00007 | 0.008 | 0.100 | 0.037 | 1.055 | 0.203 | 0.034 | 0.036 | 0.014 | 1.222 | 0.699 |
| | Minimum spacing metric | Logistic NSGA-II | 8 | 1.01392 | 1.00055 | 0.062 | 0.329 | 0.099 | 1.048 | 0.366 | 0.078 | 0.174 | 0.216 | 1.053 | 0.502 |
| PID | Minimum hypervolume indicator | Logistic NSGA-II | 9 | 1.00726 | 1.15637 | 0.991 | 0.570 | 0.376 | - | - | 0.485 | 0.269 | 0.800 | - | - |
| | Maximum diversity metric | Logistic NSGA-II | 10 | 1.01164 | 1.00137 | 0.411 | 0.375 | 0.005 | - | - | 0.316 | 0.149 | 0.042 | - | - |
| | Maximum Pareto spread | Logistic NSGA-II | 11 | 1.01164 | 1.00137 | 0.411 | 0.375 | 0.005 | - | - | 0.316 | 0.149 | 0.042 | - | - |
| | Minimum spacing metric | Logistic NSGA-II | 12 | 1.01378 | 1.00095 | 0.322 | 0.250 | 0.006 | - | - | 0.300 | 0.178 | 0.052 | - | - |

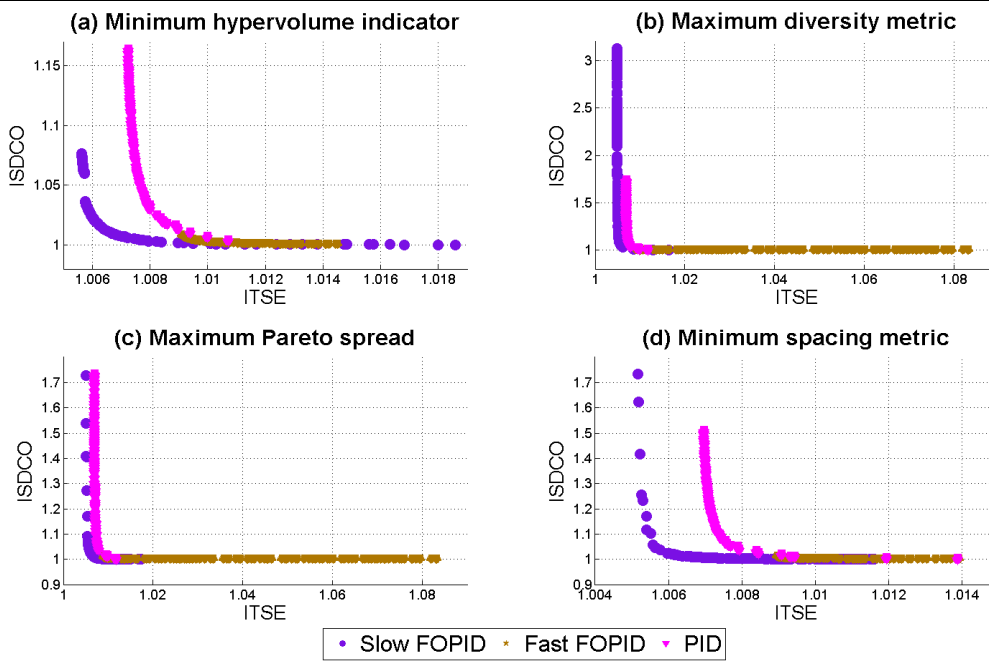

Figure 12: Comparison of non-dominance amongst the controllers using best MOO algorithms for each criterion.

## 6.2. Time domain performance of the LFC system

The time domain performance of the best compromise controller parameters reported in Table 2 are now compared for the four Pareto metrics (hypervolume indicator, diversity metric, Pareto spread and spacing metric). Since depending on the spread of the Pareto front in Figure 12 the best compromise solution may be have different time domain characteristics. The time domain responses of the grid frequency oscillation in both the areas, tie line power flow and the control signals have been reported for each Pareto metrics in Figure 13-Figure 16.

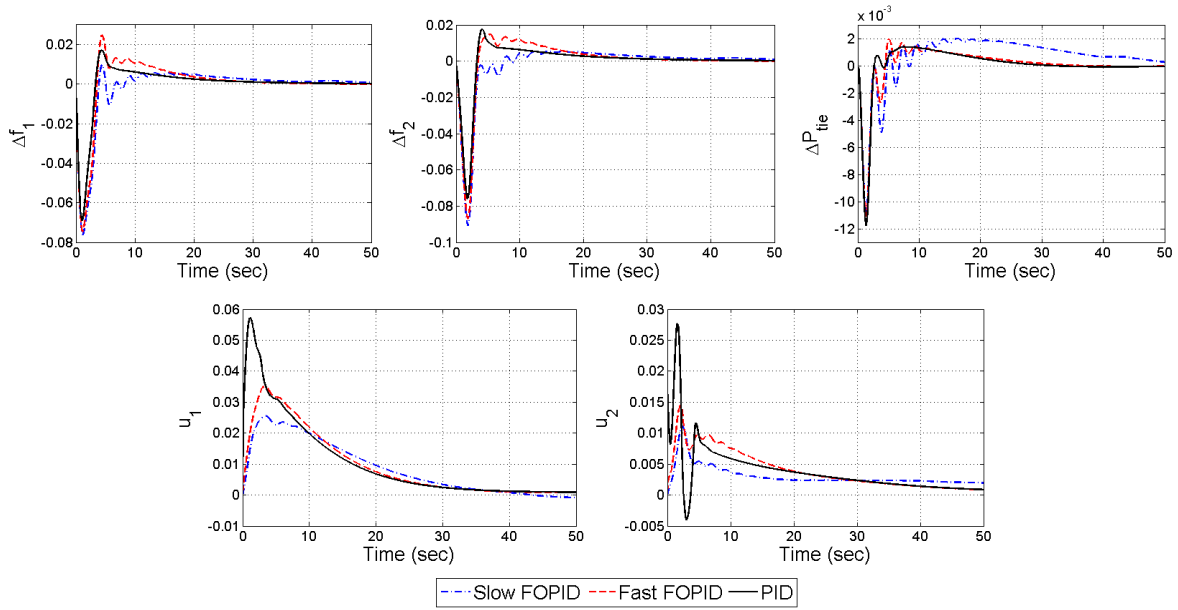

**Figure 13:** Performance comparison of the best compromise solution of the three type of nondominated controllers obtained using the minimum hypervolume indicator criterion.

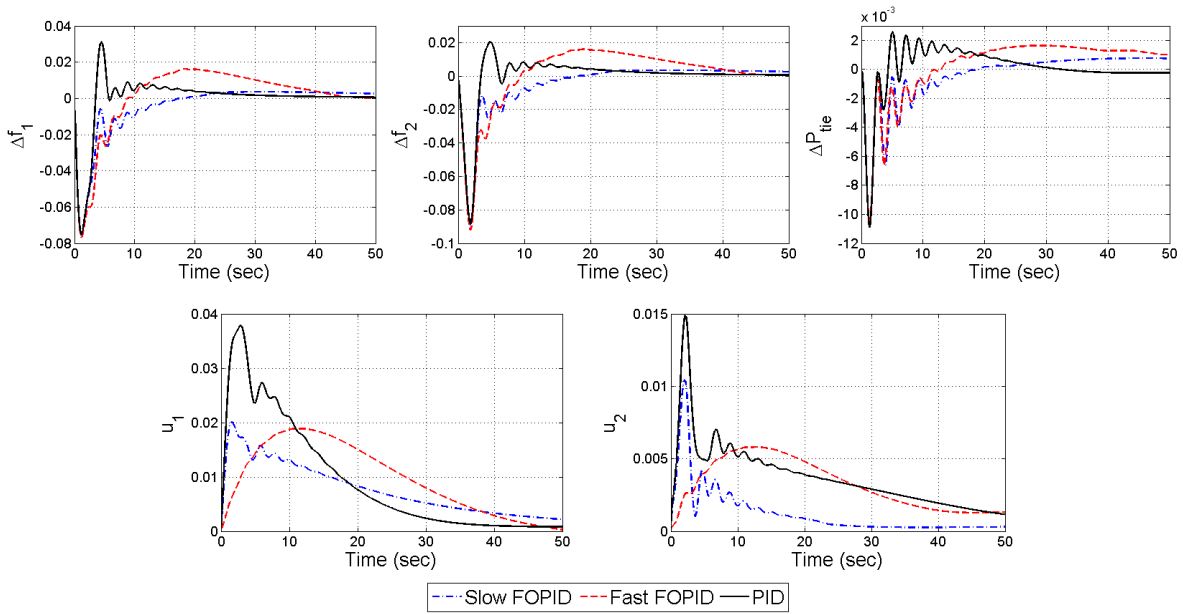

**Figure 14:** Performance comparison of the best compromise solution of the three type of nondominated controllers obtained using the maximum diversity metric criterion.

The time responses of the power system with the fuzzy based best compromise solution have been shown in Figure 13 according to the solution obtained with hypervolume indicator criterion. In other words, these solutions for the three different controller structures are obtained by calculating the best compromise solution with the fuzzy based mechanism for each of the Pareto fronts in Figure 12(a). It is observed that the time domain performance of the PID (for supressing the oscillations in $\Delta f_1, \Delta f_2$ and $\Delta P_{tie}$) lies between those of the slow and the fast FOPID, but the control signal required by the PID controller is much higher. The slow FOPID results in less oscillations and overshoot and also has a smaller control signal. Hence it outperforms the other two controller structures. The next exploration tries to understand whether a similar response is observed if the controllers are selected based on a different Pareto metric.

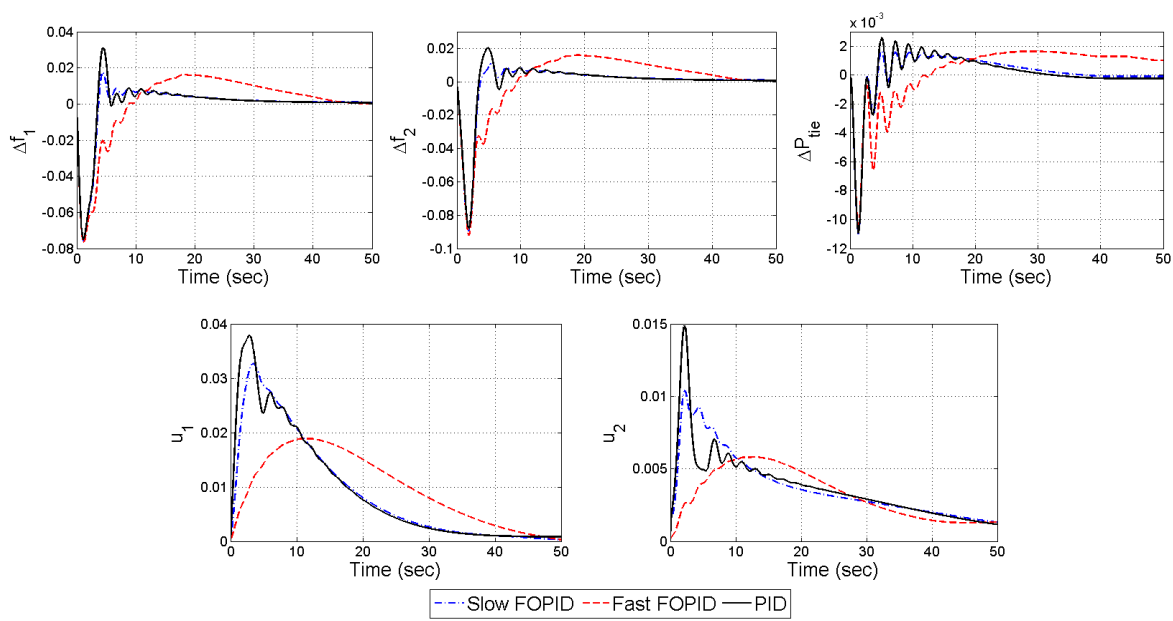

**Figure 15:** Performance comparison of the best compromise solution of the three type of nondominated controllers obtained using the maximum Pareto spread criterion.

Figure 14 indicates towards a similar conclusion which shows that the time domain performance of the fuzzy based best compromise solution obtained from the Pareto fronts of Figure 12(b), which is based on the maximum diversity metric criterion. The PID controller is found to have larger oscillations and overshoot in the time domain performance of $\Delta f_1, \Delta f_2$ and $\Delta P_{tie}$ along with a large value of control signal. The slow FOPID controller is found to outperform both the other two controller structures.

Figure 15 shows the time domain performance of the fuzzy best compromise solution obtained from the Pareto fronts of Figure 12(c), which is based on the maximum Pareto spread criterion. The fast FOPID controller has a smaller control signal but has a very sluggish time response. The slow FOPID controller has a faster time response than the fast FOPID but this comes at the cost of a higher control signal. However the slow FOPID is better than the PID in terms of both the peak overshoot and the control signal.

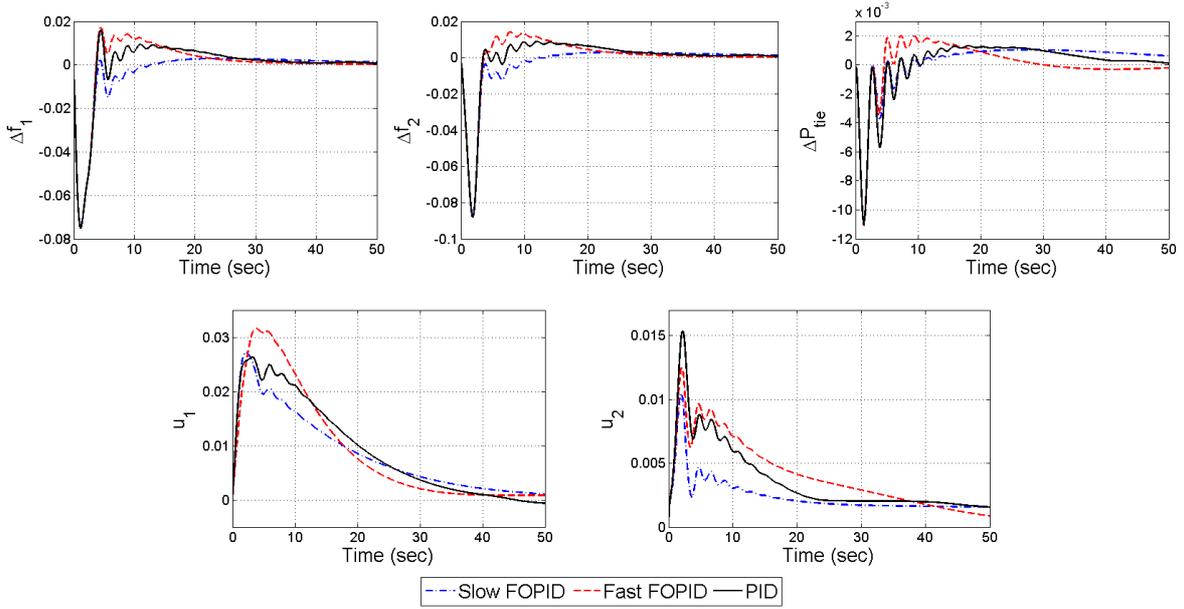

**Figure 16:** Performance comparison of the best compromise solution of the three type of nondominated controllers obtained using the minimum spacing metric criterion.

Figure 16 shows the time response of the fuzzy best compromise solution obtained from the Pareto fronts of Figure 12(d), which is based on the minimum spacing metric criterion. The slow FOPID controller is found to outperform the other two structures in terms of both time domain performance and low value of control signal.

## 7. Robustness analysis of the designed solutions

It is desirable that the designed controllers should work in a wide range of operating conditions without significant deterioration in performance. In other words, the controllers should be robust with respect to changes in system parameters. To illustrate this, the fuzzy best compromise solutions for the slow FOPID, fast FOPID and the PID controllers for two different Pareto metrics (minimum hypervolume indicator and maximum diversity metric), are simulated by varying the synchronisation coefficient ($T_{12}$). The corresponding time response curves for $\Delta f_1$, $\Delta f_2$, $\Delta u_1$, $\Delta u_2$ and $\Delta P_{tie}$ are plotted in Figure 17 and Figure 18 respectively. A similar study for the single objective FOPID controller has also been done in [34].

From the time response characteristics of $\Delta f_1$, $\Delta f_2$ and $\Delta P_{tie}$ in Figure 17 with increase in $T_{12}$ by a factor of two, it might appear that the PID controller performs better than the FOPID versions as the latter introduces small oscillations and do not settle to a steady state value quickly. However the control signals are drastically higher for the PID controller and have sharp jumps which might be detrimental for the governor.

Figure 18 shows the robustness of the fuzzy best compromise solutions for the case of the maximum diversity metric criterion with increase in $T_{12}$ by a factor of two and three. It can also be observed that the fast FOPID controller is better than the slow FOPID and the PID controller when $T_{12}$ is increased gradually. For the latter two controllers the system becomes unstable, while the fast FOPID controller is still able to maintain a time response almost similar to the nominal case which proves the superiority of the FOPID in LFC over the PID.

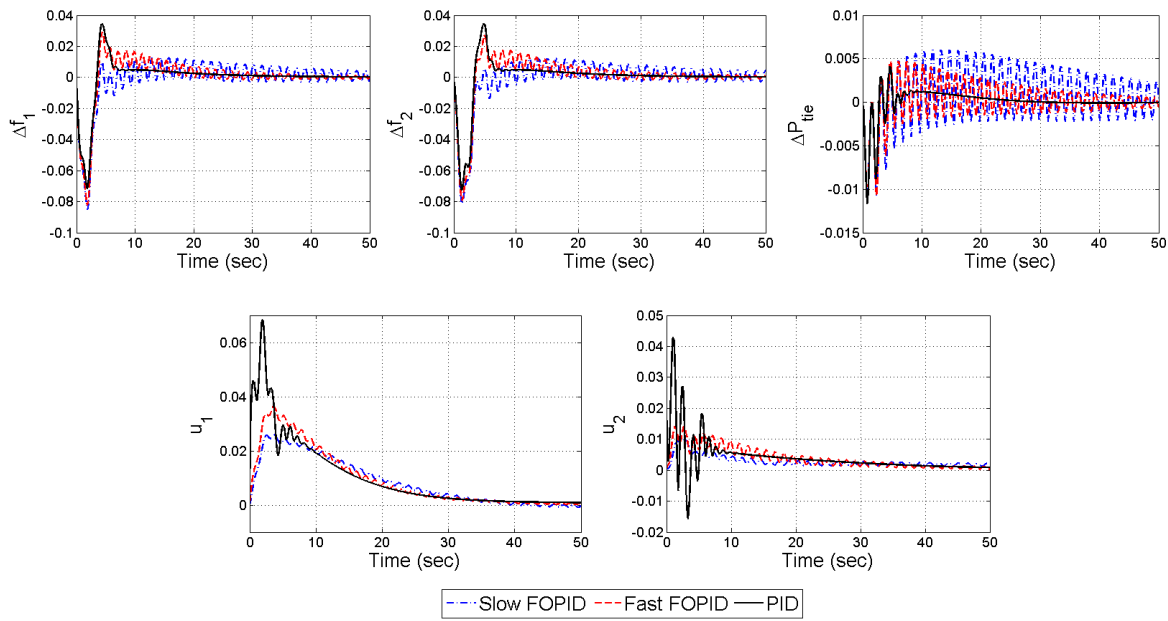

**Figure 17:** Effect of two times increase in $T_{12}$ with the best compromise solution of the minimum hypervolume indicator criterion

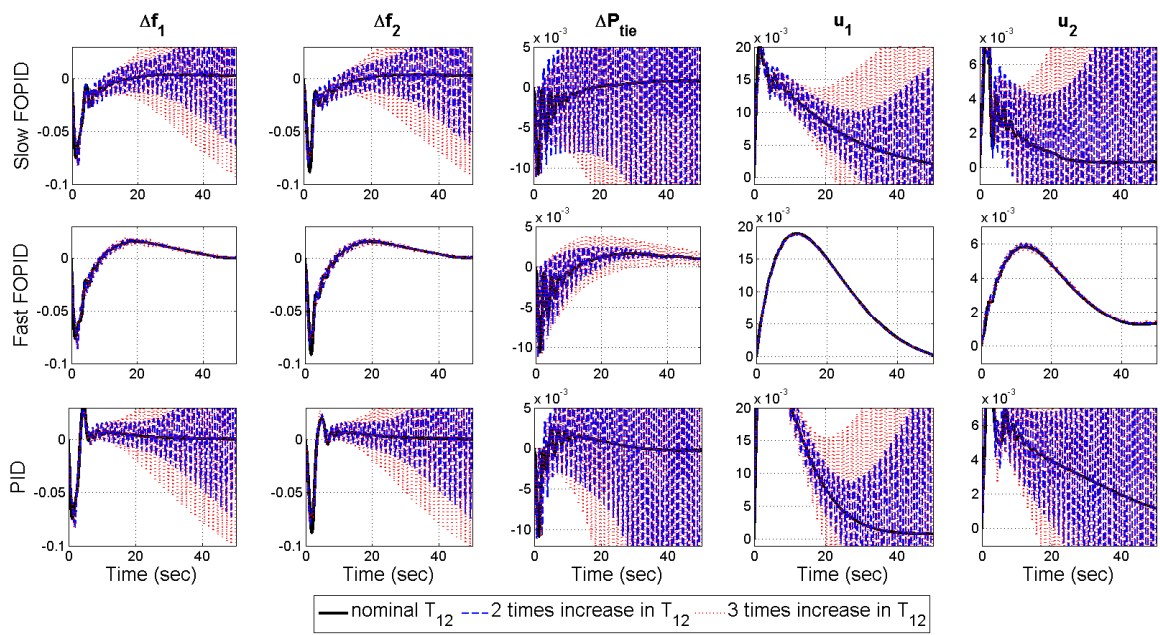

**Figure 18:** Effect of gradual increase in $T_{12}$ with the best compromise solution of the maximum diversity metric criterion.

Simulation reported so far has been done with the nominal system parameters (in section 6) and under uncertain parameters of the power system (section 7). The capability of the three controller structures to damp grid frequency oscillations even in the presence of a random change in load pattern in both the areas [29] is explored next. Figure 19 shows the time response of the system obtained by different controllers where the load patterns are randomly changing within $P_{L1} = 0.02\,pu$ and $P_{L2} = 0.008\,pu$ [37]. In terms of fast settling time and low controller effort, the slow FOPID controller is found to be better than the other two controller structures.

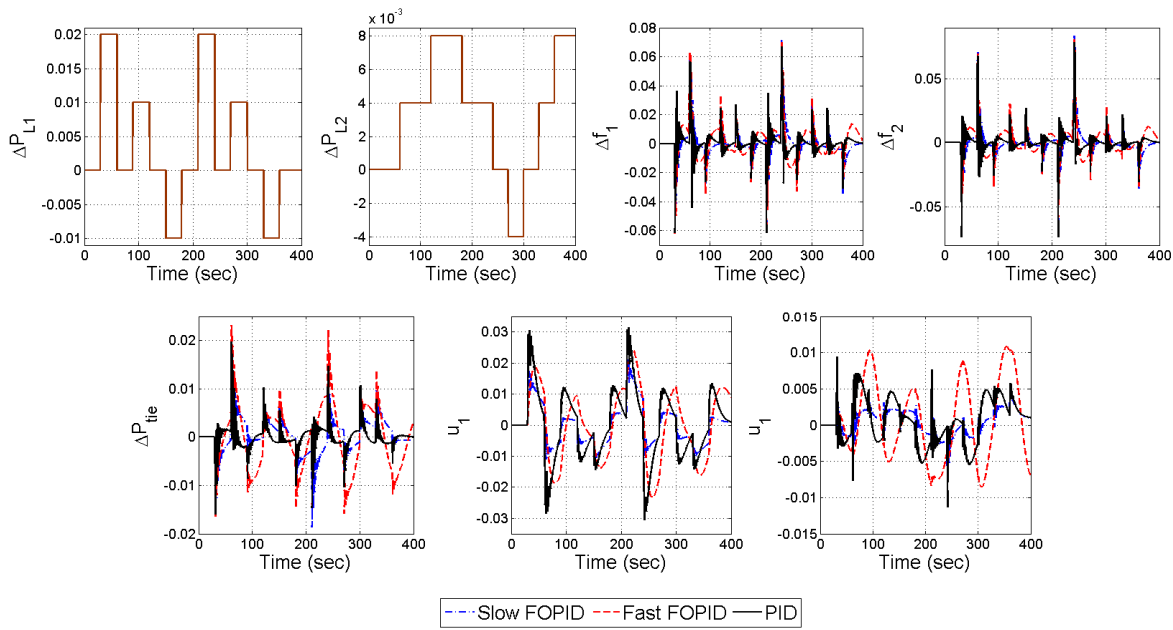

**Figure 19:** Effect of random load change in both the areas with the best compromise solution of the maximum diversity metric criterion.

Overall, both the designed solutions using the PID and the FOPID controllers show sufficient robustness to system parameter variations and random load change with FOPID variants outperforming the PID. Therefore the FOPID as load frequency controller could be applied in a practical setting where significant uncertainty exists with respect to system parameters as well as the change in load pattern.

## 8. Discussions and Conclusions

The following points summarise the main findings of the reported simulations in the paper.

- Irrespective of the chosen MOO metric (like hypervolume indicator, total Pareto spread etc.), the Pareto front obtained by the chaotic NSGA-II algorithm always results in a better set of solutions (in terms of non-domination). In other words, even though all the algorithms produce a set of non-dominated solutions as the output, those that are obtained from the chaotic NSGA II are more non-dominated vis-à-vis those obtained from their non-chaotic counterparts. For the present LFC problem, in most cases the logistic map assisted NSGA-II works better than the corresponding Henon map assisted version unlike [43].

- Under nominal conditions and random load change, slow FOPID performs better in terms of control system performance as indicated from the fast settling time and keeping the maximum values of $\Delta f_1$, $\Delta f_2$, $\Delta u_1$, $\Delta u_2$ and $\Delta P_{tie}$ low.

- When the system parameters are perturbed (e.g. synchronizing coefficient is changed), either the fast or the slow FOPID controller is better depending on the chosen Pareto metric, but it is always better than the PID controller.

In this paper, multi-objective design of an FOPID controller is done for LFC of a two area power system with GRC in turbine, reheater stages and dead-band in the governor. The NSGA-II algorithm and its chaotic versions are employed for the MOO task. Different Pareto metrics are calculated and the optimization algorithms are evaluated based on multiple Pareto metrics taken together. Numerical simulations show that the chaotic versions of the NSGA-II algorithm gives better Pareto solutions over the ordinary NSGA-II algorithm. In general, the chaotic logistic map assisted NSGA-II performs better over the chaotic Henon map assisted NSGA-II. It is also shown that the FOPID controller outperforms the PID controller for multi-objective designs of the two area LFC problem. Thus the FOPID controllers are a viable alternative to the conventional IO-PID controllers in load frequency control of inter-connected power systems.